\documentclass[11pt,draft,twoside]{article}
\usepackage[cp1251]{inputenc}
\usepackage{latexsym,amssymb,amsmath}
\title{The commutative Moufang loops with minimum conditions for
subloops II}

\author{N. I. SANDU}
\date{}
\begin{document}
\maketitle
\begin{abstract}
It is proved that the following conditions are equivalent for an
 infinite non-associative commutative Moufang loop $Q$: 1) $Q$
 satisfies the minimum condition for subloops; 2) if the loop $Q$
 contains a centrally solvable subloop of class $s$, then it
 satisfies the minimum condition for centrally solvable subloops
 of class $s$; 3) if the loop $Q$ contains a centrally nilpotent
 subloop of class $n$, then it satisfies the minimum condition for
 centrally nilpotent subloops of class $n$; 4) $Q$ satisfies the
 minimum condition for non-invatiant associative subloops. The
 structure of the commutative Moufang loops, whose infinite
 non-associative subloops are normal, is examined.
\smallskip\\
 \end{abstract}

This paper is the continuation of the article [1], where are
examined the construction of the commutative Moufang loops
(abbreviated CML) with the minimum condition for subloops. A
normal weakening for this condition is the minimum condition for
the centrally solvable (centrally nilpotent) subloops of a given
class. A broader question regarding these conditions is examined
in section 2, and namely, the existence in a CML of infinite
centrally solvable (centrally nilpotent) subloops, possessing a
property, which, by analogy with the group theory [2], will be
called steady centrally solvability (steady centrally nilpotence).
We will say that an infinite centrally solvable (centrally
nilpotent) of the class $n$ of the loop $Q$ is \textit{steadily
centrally solvable} (\textit{steadily centrally nilpotent}), if
any infinite centrally solvable (centrally nilpotent) subloop of
the class $n$ of loop $Q$ contains a proper subloop of central
solvability (central nilpotence) of class $n$. It turned out that
the existence of steadily centrally solvable (centrally nilpotent)
subloop of a certain given class $n$ in CML is equivalent to the
existence of an infinite decreasing series of subloops in CML. In
particular it follows from here that for a CML, possessing a
centrally solvable (centrally nilpotent) subloop of a certain
class $n$, the minimum condition for subloops is equivalent to the
minimum condition for subloops, which have the same class of
central solvability (central nilpotence) $n$.

It is shown in section 3 that the minimum condition for subloops
and for non-invatiant associative subloops are equivalent in the
infinite non-associative CML. The infinite non-associative CML,
which do not have proper infinite non-associative subloops are
described in section 2. A weakening of the last condition is the
condition for infinite non-associative CML, when all infinite
subloops are normal in them . The construction of such CML is
given in section 4.

\section{Preliminaries}

A \textit{multiplication group} $\frak M(Q)$ of the CML $Q$ is a
group generated by all \textit{translations} $L(x)$, where $L(x)y
= xy$. The subgroup $I(Q)$ of the group $\frak M(Q)$, generated by
all the \textit{inner mappings} $L(x,y) =L(xy)^{-1}L(x)L(y)$ is
called a \textit{inner mapping group} of the CML $Q$. The subloop
$H$ of the CML $Q$ is called \textit{normal} (\textit{invariant})
in $Q$, if $I(Q)H = H$.
\smallskip\\
\textbf{Lemma 1.1} [3]. \textit{The inner mappings are
automorphisms in the commutative Moufang loops.}

Further we will denote by $<M>$ the subloop of the loop $Q$,
generated by the set $M \subseteq Q$.
\smallskip\\
\textbf{Lemma 1.2} [3]. \textit{Let $H$ and $K$ be such loop's
subloops, that $K$ is normal in $<H, K>$. Then $HK = KH =  <H,
K>$.}

The \textit{associator} $(a,b,c)$ of the elements $a, b, c$ of the
CML $Q$ is defined by the equality $ab\cdot c = (a \cdot
bc)(a,b,c)$. The identities:

$$L(x,y)z = z(z,y,x),\eqno{(1.1)}$$

$$(x,y,z) = (y^{-1},x,z) = (y,x,z)^{-1} = (y,z,x),\eqno{(1.2)}$$

$$(x^p,y^r,z^s) = (x,y,z)^{prs}, \eqno{(1.3)} $$

$$(x,y,z)^3 = 1 \eqno{(1.4)}$$

$$(xy,u,v) = (x,u,v)((x,u,v),x,y)(y,u,v)((y,u,v).y,x)
\eqno{(1.5)}$$ hold  in a CML [3].
\smallskip\\
\textbf{Lemma 1.3} [3]. \textit{The periodic commutative Moufang
loop is locally finite}.
\smallskip\\
\textbf{Lemma 1.4} [4]. \textit{The periodic commutative Moufang
loop $Q$ decomposes into a direct product of its maximal
$p$-subloops $Q_p$, and besides, $Q_p$ belongs to the centre $Z(Q)
= \{x \in Q \vert (x,y,z) = 1 \forall y, z \in Q\}$ of CML $Q$ for
$p \neq 3$.}

We denote by $Q_i$ (respect. $Q^{(i)}$) the subloop of the CML
$Q$, generated by all associators of the form
$(x_1,x_2,\ldots,x_{2i+1})$ (respect.
$(x_1,\ldots,x_{3^i})^{(i)}$) where \break
$(x_1,\ldots,x_{2i-1},x_{2i},x_{2i+1}) =
((x_1,\ldots,x_{2i-1}),x_{2i},x_{2i+1})$ (respect.
$(x_1,\ldots,x_{3^i})^{(i)} = ((x_1,\ldots,x_{3^{
i-1}})^{(i-1)},(x_{3^{i-1}+1},\ldots,x_{2\cdot
3^{i-1}})^{(i-1)},(x_{2\cdot
3^{i-1}+1},\ldots,x_{3^i})^{(i-1)}))$. The series of normal
subloops $1 = Q_0 \subseteq Q_1 \subseteq \ldots \subseteq Q_i
\subseteq \ldots$ (respect. $1 = Q^{(o)} \subseteq Q^{(1)}
\subseteq \ldots \subseteq Q^{(i)} \subseteq \ldots $) is called
the \textit{lower central series} (respect. \textit{derived
series}) of the CML $Q$. We will also use for associator subloop
the designation $Q^{(1)} = Q'$.

The CML $Q$ is \textit{centrally nilpotent} (respect.
\textit{centrally solvable}) \textit{of class $n$} if and only if
its lower central series (respect. derived series) has the form $1
\subset Q_1 \subset \ldots \subset Q_n = Q$ (respect. $1 \subset
Q^{(1)} \subset \ldots \subset Q^{(n)} = Q$).
\smallskip\\
\textbf{Lemma 1.5 (Bruck-Slaby Theorem)} [3]. \textit{Let $n$ be a
positive integer, $n \geq 3$. Then every commutative Moufang loop
$Q$ which can be generated by $n$ elements is centrally nilpotent
of class at most $n - 1$.}

Let $M$ be a subset, $H$ be a subloop of the CML $Q$. The subloop

$$Z_H(M) = \{x \in H \vert (x,u,v) = 1 \forall u, v \in M\}$$ is
called the \textit{centralizer} of the set $M$ in the subloop $H$.
\smallskip\\
\textbf{Lemma 1.6} [1]. \textit{If $M$ is a normal subloop of the
subloop $H$ of the commutative Moufang loop $Q$ then for $a, b \in
H$ $aZ_H(M) = bZ_H(M)$ if and only if $L(a,b)(a,u,v) = (b,u,v)$
for any $u, v \in M$.}

The \textit{upper central series} of the CML $Q$ is the series

$$1 = Z_0 \subseteq Z_1 \subseteq Z_2 \subseteq \ldots \subseteq
Z_{\alpha} \subseteq \ldots$$ of the normal subloops of the CML
$Q$, satisfying the conditions: 1) $Z_{\alpha} = \sum_{\beta <
\alpha} Z_{\beta}$ for the limit ordinal and 2) $Z_{\alpha +
1}/Z{\alpha} = Z(Q/Z_{\alpha})$ for any $\alpha$.
\smallskip\\
\textbf{Lemma 1.7} [3]. \textit{The statements: 1) $x^3 \in Q$ for
any $x \in Q$; 2) the quotient loop $Q/Z(Q)$ has the exponent 3
hold  for the commutative Moufang loop $Q$.}

The CML $Q$ is called \textit{divisible} it the equation $x^n = a$
has at least one solution in $Q$, for any $n > 0$ and any element
$a \in Q$.
\smallskip\\
\textbf{Lemma 1.8} [1]. \textit{The following conditions are
equivalent for a commutative Moufang loop $D$: 1) $D$ is a
divisible loop; 2) $D$ is a direct factor for any commutative
Moufang loop that contains it.}
\smallskip\\
\textbf{Lemma 1.9} [1]. \textit{The following conditions are
equivalent for a commutative Moufang loop $Q$: 1) $Q$ satisfies
the minimum condition for subloops; 2) $Q$ is a direct product of
a finite number of quasicyclic group, lying in the centre $Z(Q)$,
and a finite loop.}

\section{Steadily centrally solvable (centrally nilpotent)
commutative Moufang loops}

\textbf{Lemma 2.1.} \textit{The infinite centrally solvable
(centrally nilpotent) commutative Moufang loop $Q$ of class $n$
contains a proper centrally solvable (centrally nilpotent) subloop
of class $n$.}
\smallskip\\
\textbf{Proof.} Let us suppose the contrary, i.e., all proper
subloops of the centrally solvable CML $Q$ of class $n$ has a
class of central solvability, less that $n$. Let us prove that in
such a case the CML is finite.

As the CML $Q$ is centrally solvable of the class $n$, there are
such elements $a_1, \ldots, a_{3^{n-1}}$ in $Q$ that
$(a_1,\ldots,a_{3^{n-1}})^{(n-1)} \neq 1$. Due to the fact that
all proper subloops of the CML $Q$ are centrally solvable of class
less than $n$, the elements $a_1, \ldots, a_{3^{n-1}}$ generate
the CML $Q$. Without violating the unity, we will suppose that all
the elements $a_1, \ldots, a_{3^{n-1}}$ are different. For
example, let the element $a_1$ have an infinite order. Then the
subloop $<a_1^4, \ldots, a_{3^{n-1}}>$ is proper in $Q$. Now, by
the identities (1.3), (1.4) we calculate

$$(a_1^4,\ldots,a_{3^{(n-1)}})^{(n-1)} =
((a_1,\ldots,a_{3^{(n-1)}})^{(n-1)})^4 =$$ $$=
(a_1,\ldots,a_{3^{(n-1)}})^{(n-1)} \neq 1.$$

We have obtained that the proper subloop $H$ is centrally solvable
of the class $n$. Contradiction. Consequently, the generators of
the CML $Q$ have a finite orders. Based on Lemma 1.3, we conclude
that the CML $Q$ is finite. Contradiction. The second case is
proved by analogy.
\smallskip\\
\textbf{Corollary 2.2.} \textit{The centrally solvable (centrally
nilpotent) commutative Moufang loop of class $n$, whose proper
subloops have a class of central solvability (central nilpotence)
less that $n$, is a finite loop.}
\smallskip\\
\textbf{Lemma 2.3.} \textit{If a non-periodic commutative Moufang
loop $Q$ contains a finite centrally solvable (centrally
nilpotent) subloop $H$ of class $n$, then it contains a steadily
centrally solvable (centrally nilpotent) subloop of class $n$.}
\smallskip\\
\textbf{Proof.} If $a$ is an element of an infinite order, then by
lemma 1.6 $a^{3^k} \in Z(Q)$, where $k = 1, 2, \ldots, Z(Q)$ is
the centre of the CML $Q$. Then the subloop $<a^{3^k}, H>$ is
steadily centrally solvable (centrally nilpotent) of the class
$n$.
\smallskip\\
\textbf{Lemma 2.4.} \textit{Let the commutative Moufang loop $Q$,
which does not satisfy the minimum condition for subloops, be
centrally solvable (centrally nilpotent) of the class $n$. Then
$Q$ possesses a proper infinite centrally solvable (centrally
nilpotent) subloop of the class $n$.}
\smallskip\\
\textbf{Proof.} Let the infinite CML $Q$ be centrally solvable of
the class $n$ and all its proper centrally solvable subloops of
the class $n$, be finite. By Lemma 2.1 there exists a finite
proper centrally solvable subloop $K$ of the class $n$ of the
order $m$ in the CML $Q$.

If $L$ is an arbitrary normal subloop of a finite index of the CML
$Q$, then by Lemma 1.2 $LK$ is an infinite centrally solvable
subloop of the class $n$ and therefore $LK = Q$. By the relation

$$Q/L =LK/L \cong K(K \cap L)$$ the index of the normal subloop
$L$ is not greater that $m$ in the CML $Q$. Then there exists in
the CML $Q$ a normal subloop $H$ of a finite index. The subloop
$H$ does not possess proper normal subloops of finite index, it
means that $H/H'$ is infinite. Therefore $H'K$ is a finite
subloop, and then the associator subloop $H'$ is also finite. Let
us show that the subloop $H$ is associative. Indeed, by Lemma 1.5
$aZ(H) \neq bZ(H)$ ($a, b \in H$) if and only if there exist such
elements $u, v \in H$ that $(a,u,v) \neq (b,u,v)$. Therefore the
centre $Z(Q)$ has a finite index in $H$. The subloop $H$ is normal
in the CML $Q$, i.e. it is invariant regarding the inner mapping
group, which consists of automorphisms (Lemma 1.1). Then it is
obvious that the subloop $Z(H)$ is normal in $Q$. We have obtained
that the CML $H$ contains a normal in $Q$ subloop of finite index.
But it contradicts the choice of subloop $H$. Consequently, $Z(H)
= H$. Further, the set $S$ of the elements of the group $H$,
having simple orders, is finite. It follows from the fact that the
subloop $<S>K$ (the subloop $<S>$ is normal in Q) is finite as a
proper centrally solvable subloop of the class $n$ of the CML $Q$.
It follows from here that $H$ is an abelian group with the minimum
condition for subgroups. The second case is proved by analogy.
\smallskip\\
\textbf{Corollary 2.5.} \textit{For an infinite centrally solvable
(centrally nilpotent) commutative Moufang loop to be steadily
centrally solvable (steadily centrally nilpotent), it is enough
that it does not contain divisible subloops different from the
unitary element.}
\smallskip\\
\textbf{Corollary 2.6.} \textit{For an infinite periodic centrally
solvable (centrally nilpotent) commutative Moufang loop $Q$ of the
class $n$ to be steadily centrally solvable (steadily centrally
nilpotent), it is necessary and sufficiently  that $Q$ does not
contain divisible subloops different from the unitary element.}
\smallskip\\
\textbf{Proof.} If $Q$ does not contain non-trivial divisible
subloops, then the necessity follows from the Corollary 2.5.
Conversely, for example, let the CML $Q$ be steadily centrally
solvable and let $H$ be the maximal divisible subloop of the CML
$Q$. By Lemma 1.7 $H \subseteq Z(Q)$. If $L$ is a finite centrally
solvable subloop of the class $n$, $K$ is a quasicyclic group from
$H$, then the subloop $<L, K>$ is centrally solvable of the class
$n$ and satisfies the minimum condition for subloops. By the
mentioned above and by Lemma 2.1 it is easy to show that there
exists an infinite centrally solvable subloop of the class $n$ in
the $<L, K>$, whose all subloop's proper centrally solvable
subloop of the class $n$ are finite. But it contradicts the fact
that $Q$ is steadily centrally solvable. The second case is proved
by analogy. This completes the proof of Corollary 2.6.

Let us remark that the request on the periodicity of the CML $Q$
in the Corollary 2.6 is essential (example: the additive group of
rational numbers).

We will call a \textit{minimal CML} of central solvability
(central nilpotence) class $n$ any centrally solvable (centrally
nilpotent) CML, whose all proper subloops have a class of central
solvability (central nilpotence) less that $n$. It follows the
Lemmas 2.1 and 1.4 that these are commutative Moufang $3$-loops.
\smallskip\\
\textbf{Corollary 2.7.} \textit{For a commutative Moufang loop $Q$
to be infinite centrally solvable (centrally nilpotent) of the
class $n$, and all its proper centrally solvable (centrally
nilpotent) subloops of the class $n$ to be finite, it is necessary
and sufficiently that the loop $Q$ is a direct product of
quasicyclic groups and the minimal CML of the central solvability
(central nilpotence) class $n$.}
\smallskip\\
\textbf{Proof.} We will examine only the case of central
solvability. If the infinite CML $Q$ is centrally solvable of
class $n$ and all its proper centrally solvable class $n$ are
finite, then by Lemma 2.4 $Q$ satisfies the minimum condition for
subloops. By Lemma 1.9 $Q$ decomposes into a direct product of
finite number of quasicyclic groups and a finite CML. Obviously,
if $K$ is a quasicyclic group and $L$ is a minimal subloop of
central solvability class $n$, then $Q = K \times L$. The inverse
is obvious.
\smallskip\\
\textbf{Lemma 2.8.} \textit{Let the commutative Moufang loop $Q$,
which does not satisfy the minimum condition for subloops, be
centrally solvable (centrally nilpotent) of the class $n$. Then
$Q$ possesses a steadily centrally solvable (steadily centrally
nilpotent) subloop of class $n$.}
\smallskip\\
\textbf{Proof.} Let $Q^{(t)}$ be the last member of derived series
(lower central series)

$$Q = Q^{(0)} \supset Q^{(1)} \supset \ldots \supset Q^{(t)}
\supset \ldots \supset Q^{(n)} = 1$$ of the CML $Q$, that does not
satisfy the minimum condition for subloops. If there are no
steadily centrally solvable (steadily centrally nilpotent)
subloops of class $n$ in the CML $Q$, then by Lemma 2.1 there
exists a finite centrally solvable (centrally nilpotent) subloop
of the class $n$ in it. We denote it by $H$.

If $Q$ is an non-periodic CML, then the statement follows from the
Lemmas 2.3.

Let now the subloop $Q^{(t)}$ has no elements of infinite order.
By (1.4) the subloop $Q^{(t + 1)}$ has the exponent three and by
the supposition it satisfies the minimum condition for subloops.
Then by Lemma 1.9    $Q^{(t + 1)}$ is finite. We denote by
$L/Q^{(t + 1)}$ the subgroup of the abelian group $Q^{(t)}/Q^{(t +
1)}$, generated by all elements of simple orders. It cannot be
finite, as the group $Q^{(t)}/Q^{(t + 1)}$, and then the CML
$Q^{(t)}$ would also satisfy the minimum condition for subloops.
We denote by $Z$ the centralizer of the normal subloop $Q^{(t +
1)}$ in the CML $L$. By Lemma 1.5, if $a, b \in L$, then $aZ \neq
bZ$ if and only if there exist such elements $u, v$ from $Q^{(t +
1)}$ that $L(a,b)(a,u,v) \neq (b,u,v)$. The subloop $Q^{(t + 1)}$
is normal in $Q$, then $(a,u,v) \in Q^{(t + 1)}$. As $Q^{(t + 1)}$
is finite, $L/Z$ is finite. So, the subloop $Z$ does not satisfy
the minimum condition for subloops. Now it follows from the
relations

$$Z/(Z \cap Q^{(t + 1)}) \cong Q^{(t + 1)}Z/Q^{(t + 1)} \subseteq
L/Q^{(t + 1)}$$ that $Z/(Z \cap Q^{(t + 1)})$ is an infinite
abelian group. The subloop $Z \cap Q^{(t + 1)}$ is contained in
the centre of the CML $Z$, then $Z$ is a centrally nilpotent CML
of the class 2. It follows from here that the associator subloop
$Z'$ is an abelian group of the exponent three. If the associator
subloop $Z'$ is infinite, then $Z'H$ is an unknown subloop (the
product $Z'H$ is a subloop by Lemma 1.2, as the normality of $Z'$
in $Q$ follows from the normality of $Z$ in $Q$). But if the
associator subloop $Z'$ is finite, then the subgroup $K/Z'$ of the
of the group $Z/Z'$, generated by all elements of simple orders,
should be infinite, as $Z$ does not satisfy the minimum condition
for subgroups. The subloop $K$ is normal in $Q$ as $Z$ is normal
in $Q$ and, obviously, $K$ contains no divisible subloops
different from the unitary element. Consequently, by the Corollary
2.6 $HK$ is a steadily centrally solvable (steadily centrally
nilpotent) subloop of the class $n$.
\smallskip\\
\textbf{Lemma 2.9.} \textit{An arbitrary centrally solvable
(centrally nilpotent) commutative Moufang loop $Q$ of class $n$,
that does not satisfy the minimum condition for subloops, that
does not satisfy the minimum condition for subloops, possesses a
steadily centrally solvable (steadily centrally nilpotent)
subloops of central solvability (central nilpotence) class $t$ for
any $t = 1, 2, \ldots, n$.}
\smallskip\\
\textbf{Proof.} Let $Q$ be a centrally nilpotent CML of class $n$
and let $a_1, a_2, \ldots, a_{2n+1}$ be such elements from $Q$
that $((a_1,\ldots,a_{2i+1}),
a_{2i+2},\ldots,a_{2n-1},a_{2n},a_{2n+1}) = 1$, but
$((a_1,\ldots,a_{2i+1}),a_{2i+2},\ldots,a_{2n-1}) \neq 1$. Then
the subloop $<(a_1,\ldots,a_{2i+1}),\break
a_{2i+2},\ldots,a_{2n+1}> = H$ is centrally nilpotent of class $n
- 1 = t$. In the case of central solvability we will examine the
$(n - i)$-th member of the derived series $Q^{(n-i)}$ instead of
$H$.

If the subloop $H$ is not steadily centrally solvable (steadily
centrally nilpotent) of class $t$, then by Lemma 2.1 the subloop
$H$ is finite. Let the CML $Q$ not be periodic. Then by Lemma 2.3
$Q$ contains a steadily centrally solvable (steadily centrally
nilpotent) subloop of class $t$.

Let us suppose that $Q$ is a periodic CML. Let $Q^{(i)}$ be the
last member of the derived series (of the lower central series)

$$Q = Q^{(0)} \supset Q^{(1)} \supset \ldots \supset Q^{(i)}
\supset \ldots \supset Q^{(n)} = 1$$ of the CML $Q$, that does not
satisfy the minimum condition for subloops. The subloop
$Q^{(i+1)}$ satisfies the minimum condition for subloops and by
(1.4) it has the index three. Then by Lemma 1.9    it is finite.
We denote by $K/Q^{(i-1)}$ the subgroup of the abelian group
$Q^{(i)}/Q^{(i+1)}$ generated by all elements of simple orders.
The group $K/Q^{(i+1)}$ is infinite, as the CML $Q^{(i)}$ does not
satisfy the minimum condition for subloops. Let us suppose that $L
= KH, L_0 = Q^{(i+1)}H$. We remind that in the case of central
solvability $Q^{(t)} = H$. But if $Q^{(i+1)}$ is a member of the
lower central series, then the subloop $L_0$ is normal in $L$.
Indeed, for that it is enough to show that if $x \in L_0, y, z \in
L$, then $(x,y,z) \in L_0$. Any element from $L_0$ has the form
$ah$, where $a \in Q^{(i+1)}, h \in H$, and any element from $L$
has the form $uh$, where $u \in Q^{(i)}, h \in H$. If $a \in
Q^{(i+1)}, u, v \in Q^{(i)}, h_1, h_2, h_3 \in H$, then by the
identity (1.5) the associator $(ah_1,uh_2,vh_3)$ may be presented
as a product of the factors of the form $(a,x,y), (h_1,h_2,h_3),
(u,x,y)$, where $x, y \in L$. As the subloop $Q^{(i+1)}$ is normal
in $Q$, $(a,x,y) \in Q^{(i+1)}$. Further, it is obvious that
$(h_1,h_2,h_3) \in H$. If $a \in Q^{(i)}$, then it follows from
the relation $Q^{(i)}/Q^{(i+1)} \subseteq Z(Q/Q^{(i+1)})$ that
$(u,x,y) \in Q^{(i+1)}$. Consequently, the subloop $L_0$ is normal
in $L$.

It has been already constructed such a series of elements of the
CML $L$

$$g_1, g_2, \ldots, g_r \eqno{(2.1)}$$ that the normal subloops
$L_i = <L_0, g_1, \ldots, g_i>$ form s strictly ascending series
$L_0 \subset L_1 \subset \ldots \subset L_r$ and for any $i = 1,
2, \ldots, r$ the element $g_i$ is bound by an associative law
with all elements of the CML $L_{i+1}$. Let us show that the
series (2.1) can be unlimitedly continued. We denote by $Z$ the
centralizer of the subloop $L_r$ in $L$. By Lemma 1.9       if $a,
b \in L$, then $aZ \neq bZ$ if and only if there exist such
elements $u, v$ from $L_r$, that $L(a,b)(a,u,v) \neq (b,u,v)$. The
CML $L_r$ is finite and normal in $L$, therefore it is easy to see
that $L/Z$ is a finite CML. Then $Z/(Z \cap L_r)$ is an infinite
CML. Let $g_{r+1} \in Z \backslash (Z \cap L_r)$. Then $L_r \neq
<L_r, g_{r+1}> = L_{r+1}$ and the element $g_{r+1}$ is bound by an
associative law with all elements of the subloop $L_r$. So, the
series (2.1) can be unlimitedly continued. The subloop $<H, g_1,
g_2, \ldots >$ is centrally solvable (centrally nilpotent) of
class $n$ and does not satisfy the minimum condition for subloops.
Indeed, according to the choice of the element $g_i$, the quotient
loop $L_0<g_1,\ldots,g_i,\ldots>/L_0$ is infinite, and therefore
it does not satisfy the minimum condition for subloops.
Consequently, the quotient loop

$$<g_1,\ldots,g_i,\ldots>/(<g_1,\ldots,g_i,\ldots>\cap L_0)$$ does
not satisfy the minimum condition for subloops as well, and as
$L_0$ is a finite CML, the subloop $<H, g_1, \ldots, g_i, \ldots
>$ does not satisfy the minimum condition for subloops. It follows
from Lemma 2.8 that there exists on $<H, g_1, \ldots, g_i, \ldots
>$ an unknown steadily centrally solvable (steadily centrally
nilpotent) subloop of class $n$.
\smallskip\\
\textbf{Corollary 2.10.} \textit{For all centrally solvable
(centrally nilpotent) of class $n$ ($n \geq 2$) subloops of the
commutative Moufang loop $Q$, that has such a subloop, to be
steadily centrally solvable (steadily centrally nilpotent) it is
enough that all its infinite centrally solvable (centrally
nilpotent) of class $n - 1$ are steadily centrally solvable
(steadily centrally nilpotent).}
\smallskip\\
\textbf{Proof.} Let $L$ be an arbitrary infinite centrally
solvable (centrally nilpotent) of class $n$ subloop of the CML
$Q$. If $L$ is not steadily centrally solvable (steadily centrally
nilpotent), then there exist in the CML $L$ an infinite centrally
solvable (centrally nilpotent) subloop $H$ of class $n$, whose all
proper subloops of central solvability (central nilpotence) class
$n$ are finite. By Lemma 2.9 the CML $H$ satisfies the minimum
condition for subloops, and by Lemma 1.9 $H = D \times K$, where
$D$ is a divisible CML, lying in the centre $Z(H)$ and $K$ is a
finite CML. The CML $K$ is centrally solvable (centrally
nilpotent) of class $n$. Then it has an proper subloop $T$ of
central solvability (central nilpotence) class $n - 1$. The
subloop $T \times D$ is infinite centrally solvable (centrally
nilpotent) subloop of class $n - 1$, satisfying the minimum
condition for subloops. It follows from Lemma 2.9 $T \times D$ is
not steadily centrally solvable (steadily centrally nilpotent).
Contradiction.
\smallskip\\
\textbf{Corollary 2.11.} \textit{For all infinite centrally
solvable (centrally nilpotent) subloops of the commutative Moufang
loop $Q$ to be steadily centrally solvable (steadily centrally
nilpotent) is necessary and sufficient that $Q$ has no quasicyclic
groups.}

The statement follows from the fact that an infinite abelian group
is steadily centrally solvable if and only if it has no
quasicyclic groups, as well as from the Corollary 2.10.
\smallskip\\
\textbf{Theorem 2.12.} \textit{If the commutative Moufang loop $Q$
possesses a centrally solvable (centrally nilpotent) subloop $S$
of class $n$ (maybe finite), then the loop $Q$ either contains a
steadily centrally solvable (steadily centrally nilpotent) subloop
of class $n$, or satisfies the minimum condition for subloops.}
\smallskip\\
\textbf{Proof.} Let us first suppose that CML $Q$ is a countable
$p$-loop and it is not centrally solvable (centrally nilpotent).
In such a case, $Q$ is the union of the countable series of finite
subloops (by Lemma 1.3 the commutative Moufang $p$-loop is locally
finite)

$$H_1 \subset H_2 \subset \ldots \subset H_k \subset \ldots,$$
where $H_i$ is a centrally solvable (centrally nilpotent) subloop
of class $n$. We denote by $L_k$ the lower layer of the centre of
the CML $H_k$. (The \textit{lower layer} of the $p$-group $G$ is
the set $\{x \in Q \vert x^p = 1\}$). Let us now examine the CML
$R = <H_1, L_2, \ldots, L_k, \ldots>$. The CML $R$ is centrally
solvable (centrally nilpotent) of class $n$. If the CML $R$ is
infinite, then is obvious that $R$ does not satisfy the minimum
condition for subloops, and by Lemma 2.9 the CML $R$ contains a
steadily centrally solvable (steadily centrally nilpotent) subloop
of class $n$. But if the CML $R$ is finite, then the CML $<L_1,
L_2, \ldots, L_k, \ldots >$ is also finite. Therefore the centre
$Z(Q)$ of the CML $Q$ is different from the unitary element. The
upper central series $Z_1 \subseteq Z_2, \subseteq \ldots
\subseteq Z_{\beta} \subseteq \ldots$ of the CML $Q$ stabilities
on a certain ordinal number $\gamma$. If $Z_{\gamma}$ is a
centrally solvable (centrally nilpotent) CML, then the CML $Q$
contains a steadily centrally solvable (steadily centrally
nilpotent) subloop of class $n$. Indeed, in this case the quotient
loop $Q/Z_{\gamma}$ is a countable $p$-loop, and is not centrally
solvable (centrally nilpotent). Then by the above-mentioned
judgements, and as the $Q/Z_{\gamma}$ is a CML without a centre,
we obtain that the CML $Q/Z_{\gamma}$ contains a steadily
centrally solvable (steadily centrally nilpotent) subloop of class
$n$. Let it be the subloop $K/Z_{\gamma}$. By the definition of
the derived's series (of the lower central series) the subloop $K$
is centrally solvable (centrally nilpotent) and it does not
satisfy the minimum condition for subloops. Then by Lemma 2.9 the
CML $K$ contains a steadily centrally solvable (steadily centrally
nilpotent) subloop of class $n$.

Let us now that $Z_{\gamma}$ is not a centrally solvable
(centrally nilpotent) subloop and let $SZ_{\alpha}$ be the first
member of the series $SZ_1 \subset SZ_2 \subset \ldots \subset
SZ_{\beta} \ldots$ which is not a centrally solvable (centrally
nilpotent) subloop. If the CML $SZ_{\beta}$ does not satisfy the
minimum condition for at least one ordinal number $\beta$ ($\beta
< \alpha$), then by Lemma 2.9 the CML $SZ_{\beta}$ contains an
unknown steadily centrally solvable (steadily centrally nilpotent)
subloop. Let us suppose now that for all $\beta$ ($\beta <
\alpha$) the subloops $SZ_{\beta}$ satisfy the minimum condition
for subloops, and let us denote by $D$ the maximal divisible
subloop of the CML $SZ_{\alpha}$. By Lemma 1.9      $SZ_{\alpha} =
D \times \overline Z_{\alpha}$, where $\overline Z_{\alpha}$ is a
reduced CML. The subloops $SZ_{\beta}$ ($\beta < \alpha$) satisfy
the minimum condition, then by Lemmas 1.8, 1.7 $SZ_{\beta} =
D_{\beta} \times \overline Z_{\beta}$, where $D_{\beta}$ are
divisible CML, $\overline Z_{\beta}$ are finite normal reduced
subloops. Consequently, $\overline Z_{\alpha}$ is the union of the
ascending series of finite normal subloops $Z_{\beta}$ ($\beta <
\alpha$). The maximal subloop $\overline M$ of the CML $\overline
Z_{\alpha}$, that has the central solvability (central nilpotence)
class $n$, cannot be finite. Indeed, it follows from the
finiteness of the subloop $\overline M$ that $\overline M \subset
\overline Z_{\beta}$ for a certain $\beta < \alpha$. We denote by
$Z$ the centralizer of the subloop $\overline Z_{\beta}$ in the
CML $\overline Z_{\alpha}$. If $a, b, \in \overline Z_{\alpha}$,
then by Lemma 1.9       $aZ \neq bZ$ if and only if the exist such
elements $u, v \in \overline Z_{\beta}$ that $L(a,b)(a,u,v) \neq
(b,u,z)$. The subloop $\overline Z_{\beta}$ is normal in $Q$ and
it is finite, therefore the centralizer $Z$ is infinite. So, there
exists an non-unitary element $w\in Z\backslash \overline M$. The
subloop $<\overline M, w>$ has the central solvability (central
nilpotence) class $n$ and is different from the subloop $\overline
M$, that contradicts the choice of $\overline M$. Thus, $\overline
M$ is an infinite CML. By the maximality of the divisible CML $D$,
the CML $\overline M$ is a steadily centrally solvable (steadily
centrally nilpotent) subloop of class $n$ by the Corollary 2.6.

Let now $Q$ be an arbitrary CML, satisfying the theorem's
conditions. If $a$ is an element of infinite order, then by Lemma
2.9 there exists in the CML $<S, a>$ a steadily centrally solvable
(steadily centrally nilpotent) subloop of class $n$.

Let $Q$ be a periodic CML, not centrally solvable (centrally
nilpotent). By Lemma 1.4 $Q$ decomposes into a direct product of
its maximal $p$-subloops $Q_p$, besides, $Q_p$ lies in the centre
of the CML $Q$ for $p \neq 3$. Then the subloop $Q_3$ is not
centrally solvable (centrally nilpotent) and such a countable
subloop can be found within it. By the above-mentioned, the latter
contains an unknown steadily centrally solvable (steadily
centrally nilpotent) subloop.
\smallskip\\
\textbf{Corollary 2.13.} \textit{The following conditions are
equivalent for a non-associative commutative Moufang loop:}

\textit{1) the loop $Q$ satisfies the minimum condition for
subloops;}

\textit{2) if the loop $Q$ contains a centrally solvable subloop
of class $s$, it satisfies the minimum condition for the centrally
solvable subloops of class $s$;}

\textit{3) if the loop $Q$ contains a centrally nilpotent subloop
of class $n$, it satisfies the minimum condition for the centrally
nilpotent subloops of class $n$;}

\textit{4) the loop $Q$ satisfies the minimum condition for the
associative subloops;}

\textit{5) the loop $Q$ satisfies the minimum condition for
non-associative subloops.}
\smallskip\\
\textbf{Corollary 2.14.} \textit{The infinite commutative Moufang
loop $Q$, possessing a centrally solvable (centrally nilpotent)
subloop $H$ of class $n$, has also an infinite subloop such type.}
\smallskip\\
\textbf{Proof.} Let $a \in Q$ be an element of infinite order. By
Lemma 1.6 $a^{3^k} \in Z(Q), k = 1, 2, \ldots$, therefore $<H,
a^{3^k}>$ is a unknown subloop. If the periodic CML $Q$ does not
satisfy the minimum condition for centrally solvable (centrally
nilpotent) subloops of class $s$, then it contains an infinite
subloop of this type, as the CML $Q$ is locally finite (Lemma
1.3). In the opposite case, by the Corollary 2.13 and Lemma 1.9
$Q = D \times K$, where $D \subseteq Z(Q), K$ is a finite CML. In
this case $D, H>$ is an unknown subloop.
\smallskip\\
\textbf{Corollary 2.15.} \textit{Any infinite commutative Moufang
loop possesses an infinite associative subloop.}

The statement follows from the Corollary 2.14 and from the fact
the CML is monoassociative.
\smallskip\\
\textbf{Corollary 2.16.} \textit{A commutative Moufang loop with
finite centrally solvable (centrally nilpotent) subloops of class
$n$, $n = 1, 2, \ldots$, is finite itself.}

The statement is equivalent to the Corollary 2.14.

In particular, the equivalence of the conditions 1), 5) of the
Corollary 2.13 means that each infinite non-associative CML has an
infinite non-associative subloop different from itself the
exception of the case when it satisfies the minimum condition for
subloops. It is clear that not any infinite CML with the minimum
condition is an exception here. It holds true indeed.
\smallskip\\
\textbf{Proposition 2.17.} \textit{The infinite non-associative
commutative Moufang loop $Q$ does not contain its proper infinite
non-associative subloops if and if it decomposes into a direct
product of quasicyclic groups, contained in the centre $Z(Q)$ of
the loop $Q$, and a finite non-associative loop, generated by
three elements.}
\smallskip\\
\textbf{Proof.} By the Corollary 2.13 the CML $Q$ satisfies the
minimum condition for subloops, then by Lemma 1.9    $Q = D \times
H$, where $D$ is a direct product of a finite number of
quasicyclic groups, $D \subseteq Z(Q)$, $H$ is a finite CML. By
the supposition about the CML $Q$, the group $D$ contains only one
quasicyclic group.

Obviously $H$ is an non-associative CML. If $H_1$ is an arbitrary
proper subloop of the CML $H$, then by Lemma 1.2 the product
$DH_1$ is an proper infinite subloop of the CML $Q$. But then
$DH_1$ and $H_1$ are associative subloops. Consequently, all
proper subloops of the CML $Q$ are associative, and it follows
from the Lemma 1.5 [3] that $H$ is generated by tree elements. Let
now the CML $Q$ have a decomposition $Q = D \times H$, possessing
these qualities, and $L$ be an arbitrary proper subloop of the CML
$Q$. Obviously $D \subseteq L$. Then it follows from the
decomposition $Q = D \times H$ that $L = D(L \cap H)$. As $L \neq
Q$, then $L \cap H \neq H$. Then the subloop $L \cap H$, as an
proper subloop of the CML $H$, is associative. Therefore it
follows from the decomposition $L = D(L \cap H)$ that the subloop
$H$ is associative.

\section{Infinite non-associative commutative Moufang loops
with minimum condition for non-invatiant associative subloops}

\textbf{Lemma 3.1.} \textit{If the element $a$ of an infinite
order or of order three of the commutative Moufang loop $Q$
generates a normal subloop, then it belongs to the centre $Z(Q)$
of loop $Q$.}
\smallskip\\
\textbf{Proof.} If the element $1 \neq a \in Q$ generates a normal
subloop, then $L(u,v)a = a^k$ for a certain natural number $k$ and
for arbitrary fixed elements $u, v \in Q$. By (1.1) $a(a,v,u) =
a^k, (a,v,u) = a^{k-1}$. If $k = 1$, then $(a,v,u) = 1$. Therefore
$a \in Z(Q)$. Let us now suppose that $k
> 1$. Let $a^3 = 1$. Then $k = 2$ and  by (1.5) and Lemma
1.5 $a = (a,v,u), a = ((a,v,u),v,u) = 1$. We have obtained a
contradiction, as $a \neq 1$. But if $a$ has an infinite order,
then by (1.4) $(a^{k-1})^3 = (a,v,u)^3 = 1$. We have obtained a
contradiction again. Therefore the case of $k > 1$ is impossible.
This completes the proof of Lemma 3.1.
\smallskip\\
\textbf{Lemma 3.2.} \textit{The commutative Moufang loop $Q$,
containing an element of an infinite order is associative if and
only if the subloop, generated by any element of an infinite
order, is normal in $Q$.}
\smallskip\\
\textbf{Proof.} By Lemma 3.1 any element $a$ of an infinite order
of the CML $Q$ belongs to the centre $Z(Q)$. Let $b$ be an element
of a finite order of the CML $Q$. Obviously the product $ab$ has
an infinite order. Again by Lemma 3.1 $ab \in Z(Q)$. Further, by
(1.5) and (1.4) we have $1 = (ab,u,v) = L(a,b)(a,u,v)\cdot
L(b,a)(b,u,v) = (b,L(b,a)u,L(b,a)v)$, for $u, v \in Q$.
Consequently, $b \in Z(Q)$, but then the CML $Q$ is associative.
\smallskip\\
\textbf{Theorem 3.3.} \textit{If in the infinite commutative
Moufang loop $Q$ the infinite associative subloops are normal in
$Q$, then $Q$ is associative.}
\smallskip\\
\textbf{Proof.} It follows from Lemma 3.2 that it is sufficient to
examine the case when the CML $Q$ is periodic, and by Lemma 1.4 it
is sufficient to examine the case when $Q$ is a $3$-loop.

Let us now first examine the case when the CML $Q$ does not
satisfy the minimum condition for subloops. By the Corollary 4.5
from [1] none of its maximal elementary associative subloops $H$
can be finite. Let

$$H = H_1 \times H_2 \times \ldots \times H_n \times \ldots$$ be
the decomposition of the group $H$ into a direct product of cyclic
groups of order three. We denote by $Z_Q(H)$ the centralizer of
the subloop $H$ in $Q$. It is obvious that for any element $a$
from $Z_Q(H)$ there is such an infinite subgroup $H(a) \subseteq
H$ that $<a> \cap H(a) = 1$. Let $H(a) = H_1(a) \times H_2(a)$ be
a decomposition of the group $H(a)$ into a direct product of
infinite factors. As the cyclic group $<a>$ is the intersection of
the infinite associative subloops $<a>H_1(a)$ and $<a>H_2(a)$,
then $<a>$ is normal in $Q$. As the element $a$ from $Z_Q(H)$ is
arbitrary, we obtain that any subloop from $Z_Q(Q)$ is normal in
$Q$, i.e. $Z_Q(H)$ is a hamiltonian CML. Then by [4] it is an
associative subloop. Obviously, $H \subseteq Z_Q(H)$ and, as $H_i$
are cyclic groups of order three, then by lemma 3.1 $H_i \subseteq
Z(Q)$, where $Z(Q)$ is the centre of the CML $Q$. Consequently,
$Z(H) = Q$ is an associative CML.

If a CML $Q$ satisfies the minimum condition for subloops, then by
Lemma 1.9 its centre $Z(Q)$ is infinite. If $a$ is an arbitrary
element from $Q$, then the subloop $<a>Z(Q)$ is infinite and
associative. From here and from the theorem's supposition we
obtain that the subloop $<a>$ is normal in $Q$. Then by [4] the
CML $Q$ is associative.
\smallskip\\
\textbf{Lemma 3.4.} \textit{The non-periodic commutative Moufang
loop, satisfying the minimum condition for the non-invatiant cycle
groups, is associative.}
\smallskip\\
\textbf{Proof.} By Lemma 3.2 we suppose that the element $a$ of an
infinite order of the CML $Q$ generates a non-invatiant subloop.
It follows from the condition of lemma that the series

$$<a> \supset <a_t> \supset <a^{t^2}> \supset \ldots \supset
<a^{t^n}> \supset \ldots$$ should contain a normal subloop
$<a^{t^n}>$ for any natural number $t$. Let $t$ and $p$ be two
different simple numbers, $<a^{t^n}>$ and $<a^{p^k}>$ be two
normal subloops answering to them, of such a type that $u, v$ are
such integer numbers that $ut^n +vp^k = 1$. Then

$$a = a^{ut^n + vp^k} = a^{ut^n}\cdot a^{vp^k}.$$ If $x$ and $y$
are arbitrary elements from $Q$, then by Lemma 1.1 the inner
mapping $L(x,y)$ is an automorphism. Then, by the normality of the
subloops $<a^{t^n}>, <a^{p^k}>$, we obtain $L(x,y)a =
L(x,y)a^{ut^n}\cdot L(x,y)a^{vp^k} = (L(x,y)a^{t^n})^u \cdot
(L(x,y)a^{p^k})^v \in <a>$. Consequently, the subloop $<a>$ is
normal in $Q$. Contradiction. Then the CML $Q$ is associative.
\smallskip\\
\textbf{Theorem 3.5.} \textit{In non-associative commutative
Moufang loop the minimum condition for subloops and the minimum
condition for non-invatiant associative subloops are equivalent.}
\smallskip\\
\textbf{Proof.} Let us suppose that the CML $Q$, satisfying the
minimum condition for non-invatiant associative subloops, does not
satisfy the minimum condition for subloops. Then by the Corollary
2.13 the CML $Q$ does not satisfy the minimum condition for
associative subloops. Let us show that in this case the CML $Q$ is
associative, i.e. we will obtain a contradiction. By Lemma 3.4 it
is sufficient to examine the case when the CML $Q$ is periodic,
and by Lemma 1.4 when $Q$ is a $3$-loop.

As the CML $Q$ does not satisfy the minimum condition for
associative subloops, then by the Corollary 4.5 from [1] $Q$
contains an infinite direct product

$$H = H_1 \times H_2 \times \ldots \times H_n \times \ldots$$ of
cyclic groups of order three. If $a$ is an arbitrary element from
the centralizer $Z_Q(H)$ of the subloop $H$ in the CML $Q$, then
there exists such a number $n = n(a)$, that

$$<a> \cap (H_{n+1} \times H_{n+2} \times \ldots) = 1.$$ As the
CML $Q$ satisfies the minimum condition for non-invatiant
associative sub\-loops, then the infinitely descending series of
associative subloops

$$S^k(a) \supset S^{k+1}(a) \subset \ldots$$ contains a normal
subloop $S^l(a) (l = l(a))$, beginning with any natural $k \geq
n$, where $S^k(a) = <a>(H_{k+1} \times H_{k+2} \times \ldots)$. As
the intersection of all such normal subloops coincides with the
subloop $<a>$, then the latter is normal in $Q$. But $a$ is an
arbitrary element from the centralizer $Z_Q(H)$, and it means that
any normal subloop from $Z_Q(H)$ is normal. Then by [4] the CML
$Z(H)$ is associative. Further, the subgroups $H_i$ have the order
three. Then it follows from Lemma 3.1 that they belong to the
centre $Z(Q)$ of the CML Q. Then it follows from the definition of
the centralizer $Z_Q(H)$ that $Z(Q) = Q$. Consequently, the CML
$Q$ is associative.

\section{Infinite non-associative commutative Moufang loops,
in which all infinite non-associative subloops are normal}

\textbf{Lemma 4.1.} \textit{Let all infinite non-associative
subloops be normal in the infinite non-associative commutative
Moufang loop $Q$. If $H$ is an infinite non-associative subloop,
then the quotient loop $Q/H$ is a group.}
\smallskip\\
\textbf{Proof.} It is obvious that any subloop of the CML $Q$
containing $H$, is normal in $Q$. Then the quotient loop $Q/H$ is
hamiltonian, consequently by [4] it is a group.
\smallskip\\
\textbf{Proposition 4.2.} \textit{The commutative Moufang loop, in
which all its infinite non-associative subloops are normal, has a
finite associator subloop $Q'$.}
\smallskip\\
\textbf{Proof.} Let us suppose the contrary, i.e., that the
associator subloop $Q'$ is infinite. First we examine the case
when $Q'$ is non-associative. Let $H$ be a proper infinite
non-associative subloop of the CML $Q'$. Then by Lemma 4.1 $Q/H$
is a group, i.e. $Q' \subseteq H$. Contradiction. Consequently,
the associator subloop $Q'$ does not have its proper infinite
non-associative subloops. In this case, by the Corollary 2.13 the
CML $Q'$ satisfies the minimum condition for subloops. But by
(1.4) the associator subloop $Q'$ has the exponent three,
therefore it is finite.

Let us now examine the case when the infinite associator subloop
$Q'$ of the periodic CML is associative. Let $H$ be a finite
non-associative subloop of the CML $Q$. We will examine the
subloop $Q'H = \cup x_iQ', x_i \in H, i = 1, \ldots, m$. If the
infinite non-associative subloop $Q'H$ does not contain its proper
infinite non-associative subloops, then by the Corollary 2.13 it
satisfies the minimum condition for subloops. Taking into account
(1.4), it is easy to see that the CML $Q'H$ has a finite index.
Then it is finite, therefore the CML $Q'$ is also finite. It
contradicts the fact the CML $Q'H$ does not contain its proper
infinite non-associative subloops. Let $(Q'H)_1$ be the proper
infinite non-associative subloops of the CML $Q'H$. By Lemma 4.1
$Q' \subseteq (Q'H)_1$. Then $(Q'H)_1 = \cup x_iQ', i = 1, \ldots,
n, n < m$. If the infinite non-associative subloop $(Q'H)_1$ does
not contain its proper infinite non-associative subloops, then
$(Q'H)_1$ is finite, as it is shown above. Contradiction.
Therefore let $(Q'H)_2$ be the proper infinite non-associative
subloop of the CML $(Q'H)_1$. By lemma 4.1 $Q' \subseteq (Q'H)_2$,
therefore $(Q'H)_2 \subseteq \cup x_iQ', x_i \in H, i = 1, \ldots,
r, r < n$. Applying the previous judgements to the CML $(Q'H)_2$,
after a finite number of steps we will come to an infinite
non-associative subloops $(Q'H)_i$ without proper infinite
non-associative subloops. But it contradicts the statement from
the previous paragraph. Consequently, the associator subloop $Q'$
of the CML $Q$ cannot be infinite.

Finally, let us examine the case when the CML $Q$ is non-periodic.
Obviously, the subloop $H$ of the CML $Q$ is non-associative if
and only if the subloop $HZ(Q)$ is non-associative, where $Z(Q)$
is the centre of the CML $Q$. If the infinite non-associative
subloops of the CML $Q$ are normal, then the infinite
non-associative subloops of the CML $Q/Z(Q)$ are normal as well.
By Lemma 1.9 the CML $Q/Z(Q)$ has index three, then, according to
the previous case, its associator subloop $(Q/Z(Q))'$ is finite.
If $a \in Z(Q)$, then $(au,v,w) = (u,v,w)$, for any $u, v, w \in
Q$. It is easy to see from here that the associator subloop $Q'$
is finite.
\smallskip\\
\textbf{Corollary 4.3.} \textit{If in the non-periodic commutative
Moufang loop $Q$ all the infinite non-associative subloops are
normal in $Q$, then its associator subloop is a finite associative
subloop.}
\smallskip\\
\textbf{Proof.} Let us suppose that the finite associator subloop
$Q'$ is non-associative. Let $H$ be one of its minimal
non-associative subloops, and $a$ be an element of infinite order
from $Q$. By Lemma 1.9       $a^3$ belongs to the centre of the
CML $Q$. Then by Lemma 1.2, $H<a^3>$ is an infinite
non-associative subloop. By Lemma 4.1 $Q' \subseteq H<a^3>$, and
it is impossible, if $H \neq Q'$. According to the minimality of
the non-associative CML $H$, it can be presented in the form of
the product of the normal associative subloop $L$ and the cyclic
group $<b>$. Indeed, by the Moufang theorem [3] the CML $H$ is
generated by three elements $u, v, b$. By Lemma 1.5 $Q' \neq H$.
Then $L = <Q', u, v>$ is a normal associative subloop and $H =
L\cdot <b>$. Now let us take the CML $B\cdot<a^3b>$. It is an
infinite non-associative subloop and, obviously, it does not
contain $Q'$. However, by Lemma 4.1 $Q' \subseteq B$.
Contradiction. Consequently, the associator subloop $Q'$ of the
CML $Q$ is associative.
\smallskip\\
\textbf{Theorem 4.4.} \textit{If all infinite non-associative
subloops of the commutative Moufang loop $Q$ are normal in it,
then all non-associative subloops are also normal in it.}
\smallskip\\
\textbf{Proof.} Let $Q$ be an non-periodic CML and $a$ be an
element of an infinite order from $Q$. By Lemma 1.9       $a^3$
belongs to the centre of the CML $Q$. If $H$ is a finite
non-associative subloop, then by Lemma 1.2 $<a^3>H$ is an infinite
non-associative subloop from $Q$ and, consequently, it is normal
in $Q$. Therefore, $H$ is normal in $Q$.

Let now $Q$ be a periodic CML and let us suppose that the finite
non-associative subloop $L$ is not normal in $Q$. The associator
subloop $Q'$ is a normal subloop in $Q$. Therefore, by Lemma 1.9
the centralizer $Z_Q(H)$ of the subloop $H$ in $Q$ will be normal
subloop in $Q$. Let us examine the set

$$C(H) = \{x \in Z_Q(H) \vert (x,u,v) = 1 \forall u \in Z_Q(H),
\forall v \in H\}.$$ Using the identity (1.5), it is easy to show
that $C(H)$ is a subloop. Moreover, it follows from the normality
of the subloops $H, Z_Q(H)$, and by Lemma 1.1, that $C(H)$ is
normal in $Q$. Indeed, if $xC(H) = yC(H)$, then $xy^{-1} \in C(H),
(xy^{-1},u,v) = 1$ for all $u \in Z_Q(H), v \in H$. Now we will
use the identities (1.5), (1.1) and (1.3). We have $1 =
(xy^{-1},u,v) = L(x,y^{-1})(x,u,v)\cdot L(y^{-1},x)(y^{-1},u,v) =
(x, L(x,y^{-1})u, L(x,\break y^{-1})v)(y^{-1}, L(y^{-1},x)u,
L(y^{-1},x)v) \equiv (x,\overline u,\overline v)(y^{-1},\overline
u,\overline v) = (x,\overline u,\overline v)(y,\overline
u,\overline v)^{-1},\break (x,\overline u,\overline v) =
(y,\overline u,\overline v)$ for all $u \in Z_Q(H), v \in H$. It
can be proved by analogy that it follows from the equality
$(x,\overline u,\overline v) = (y,\overline u,\overline v)$ from
all $u \in Z_Q(H), v \in H$ that $xC(H) = yC(H)$. By the
Proposition 4.2 the associator subloop $Q'$ is finite. Then the
normal subloop $C(H)$ has a finite index in $Q$.

Let us show that the CML $Q$ satisfies the minimum condition for
subloops. Let us suppose the contrary. Then the subloop $C(H)$,
possessing a finite index in $Q$, does not satisfy this condition
as well. Therefore, the CML $C(H)$ has an infinite associative
subloop $K$, which decomposes into a direct product of cyclic
groups of simple orders. Otherwise, by the Corollary 2.13 and
regardless the supposition, the CML $Q$ would satisfy the minimum
condition for subloops. It is obvious that there can be emphasized
an infinite subgroup $R$, that intersects with $L$ on the unitary
element. Let $R = R_1 \times R_2$ be the decomposition of $R$ into
a direct product of two infinite subgroups $R_1, R_2$.

If $S$ is an arbitrary associative subloop of he CML $C(H)$, then
the product $SL$ is a subloop. Indeed, by Lemma 1.2, the subloop
$S$ is normal in the CML $<S, L>$. The CML $<S, L>$ consists of
all ''words'', composed of the elements of the set $S \cup L$. A
word of the length 1 is an element of the set $S \cup L$.  If $u,
v$ are words of length $m, n$ respectively, then
$u^{\epsilon_1}v^{\epsilon_2}$, where $\epsilon_1, \epsilon_2 =
\pm 1$, is a word of length $\leq m + n$. It follows from the
definition of the subloop $C(H)$ that if 1) $a \in S, u \in L$; 2)
$a, u \in S, v \in L$, then $(a,u,v) = 1$. If $a \in S, u, v \in
<S,L>$ then, using (1.2), (1.5) and the associativity of the
subloop $S$, it can be proved induction on the sum of the length
of the words $u, v$ that $(a,u,v) = 1$. Then by (1.1) $L(v,u)a =
a$, i.e. the subloop $S$ is normal in $<S, L>$. Therefore $<S, L>
= SL$.

By the above prove fact, the products $R_1L, R_2L$ are subloops.
As they are infinite and non-associative, they are normal in the
CML $Q$. Then their intersection $L$ is also a normal subloop in
$Q$. We have obtained a contradiction despite the supposition of
the noninvariance  of the subloop $L$. In this case, by Lemma 1.8
the CML $Q$ decomposes into a direct product of the divisible
group $D$, lying in the centre $Z(Q)$ of the CML $Q$ and the
finite CML $M$. If $L \neq M$, then the product $DL$ is an
infinite non-associative subloop of the CML $Q$, therefore the
subloop $L$ is also normal in $Q$. We have obtained a
contradiction of the fact that $L$ is not normal in $Q$. This
completes the proof of Theorem 4.4.

By the Corollary 4.3 a non-periodic CML, whose infinite
non-associative sub\-loops, are normal in it, has a finite
associative associator subloop. The follows statement holds true
for the general case.
\smallskip\\
\textbf{Corollary 4.5.} \textit{If all (infinite) non-associative
subloops of the (infinite) non-associative commutative Moufang
loop $Q$ are normal in it, then its associator subloop $Q'$ is
centrally nilpotent, and the loop $Q$ itself is centrally solvable
of a class not greater than three.}
\smallskip\\
\textbf{Proof.} By the Proposition 4.2, the associator subloop
$Q'$ is finite. Then by Lemma 1.5 $Q'$ is centrally nilpotent.

Let us suppose that the second associator subloop $Q^{(2)}$ of the
CML $Q$ is non-associative. Then any subloop that contains
$Q^{(2)}$ is non-associative, and by Theorem 4.4, it is normal in
$Q$. Obviously, the CML $Q/Q^{(2)}$ is hamiltonian, when it is an
abelian group, by [4]. Therefore, $Q' \subseteq Q^{(2)}$, i.e. $Q'
= Q^{(2)}$. But the associator subloop $Q'$ is centrally
nilpotent, therefore $Q' \neq Q^{(2)}$. Contradiction.
Consequently, $Q^{(2)}$ is an associative subloop, and the CML $Q$
is centrally solvable of step not greater than three.

\smallskip

Tiraspol State University, Moldova

str. Iablochkin 5, Chishinau, MD-2069

Moldova

E-mail: sandumn@yahoo.com
\end{document}